\documentclass{amsart}
\usepackage{amsfonts}

\newtheorem{theorem}{Theorem}
\theoremstyle{plain}

\newtheorem{corollary}{Corollary}

\newtheorem{lemma}{Lemma}

\newtheorem{proposition}{Proposition}
\newtheorem{remark}{Remark}

\numberwithin{equation}{section}
\input{tcilatex}

\begin{document}
\title[Gr\"{u}ss Type Inequalities]{Gr\"{u}ss Type Discrete Inequalities in
Normed Linear Spaces, Revisited}
\author{Sever S. Dragomir}
\address{School of Computer Science and Mathematics\\
Victoria University of Technology\\
PO Box 14428, MCMC 8001\\
Victoria, Australia.}
\email{sever@matilda.vu.edu.au}
\urladdr{http://rgmia.vu.edu.au/SSDragomirWeb.html}
\date{22 May, 2003.}

\begin{abstract}
Some sharp inequalities of Gr\"{u}ss type for sequences of vectors in real
or complex normed linear spaces are obtained. Applications for the discrete
Fourier and Mellin transform are given. Estimates for polynomials with
coefficients in normed spaces are provided as well.
\end{abstract}

\maketitle

\section{Introduction}

The following Gr\"{u}ss type inequalities for vectors in normed linear
spaces are known.

\begin{theorem}
\label{t1.1}Let $\left( X,\left\Vert \cdot \right\Vert \right) $ be a normed
linear space over the real or complex number field $\mathbb{K}$ $\left( 
\mathbb{K}=\mathbb{R},\mathbb{C}\right) ,$ $\overline{\mathbf{\alpha }}%
=\left( \alpha _{1},\dots ,\alpha _{n}\right) \in \mathbb{R}^{n},$ $%
\overline{\mathbf{p}}=\left( p_{1},\dots ,p_{n}\right) \in \mathbb{R}_{+}^{n}
$ with $\sum_{i=1}^{n}p_{i}=1$ and $\overline{\mathbf{x}}=\left( x_{1},\dots
,x_{n}\right) \in X^{n}.$ Then one has the inequalities:%
\begin{align}
0& \leq \left\Vert \sum_{i=1}^{n}p_{i}\alpha
_{i}x_{i}-\sum_{i=1}^{n}p_{i}\alpha _{i}\cdot
\sum_{i=1}^{n}p_{i}x_{i}\right\Vert   \label{1.1} \\
& \leq \left\{ 
\begin{array}{l}
\left[ \sum\limits_{i=1}^{n}i^{2}p_{i}-\left(
\sum\limits_{i=1}^{n}ip_{i}\right) ^{2}\right] \max\limits_{1\leq j\leq
n-1}\left\vert \Delta \alpha _{j}\right\vert \max\limits_{1\leq j\leq
n-1}\left\Vert \Delta x_{j}\right\Vert ,\text{ \ \cite{DB}}; \\ 
\\ 
\dfrac{1}{2}\sum\limits_{i=1}^{n}p_{i}\left( 1-p_{i}\right)
\sum\limits_{j=1}^{n-1}\left\vert \Delta \alpha _{j}\right\vert
\sum\limits_{j=1}^{n-1}\left\Vert \Delta x_{j}\right\Vert ,\text{ \ \ \ \ \ 
\cite{D1}}; \\ 
\\ 
\sum\limits_{1\leq j<i\leq n}p_{i}p_{j}\left( i-j\right) \left(
\sum\limits_{k=1}^{n-1}\left\vert \Delta \alpha _{j}\right\vert ^{p}\right)
^{\frac{1}{p}}\left( \sum\limits_{k=1}^{n-1}\left\Vert \Delta
x_{k}\right\Vert ^{q}\right) ^{\frac{1}{q}} \\ 
\hfill \text{if }p>1,\ \frac{1}{p}+\frac{1}{q}=1,\text{ \ \ \ \cite{D2}.}%
\end{array}%
\right.   \notag
\end{align}%
The constant $1$ in the first branch, $\frac{1}{2}$ in the second branch and 
$1$ in the third branch are best possible in the sense that they cannot be
replaced by smaller constants.
\end{theorem}

The following corollary providing some inequalities for unweighted means
holds as well.

\begin{corollary}
\label{c1.2}Let $\left( X,\left\Vert \cdot \right\Vert \right) ,$ $\overline{%
\mathbf{\alpha }}\in \mathbb{K}^{n}$ and $\overline{\mathbf{x}}\in X^{n}$ be
as in Theorem \ref{t1.1}. Then one has the inequalities%
\begin{align}
0& \leq \left\Vert \frac{1}{n}\sum_{i=1}^{n}\alpha _{i}x_{i}-\frac{1}{n}%
\sum_{i=1}^{n}\alpha _{i}\cdot \frac{1}{n}\sum_{i=1}^{n}x_{i}\right\Vert 
\label{1.2} \\
& \leq \left\{ 
\begin{array}{l}
\dfrac{1}{12}\left( n^{2}-1\right) \max\limits_{1\leq j\leq n-1}\left\vert
\Delta \alpha _{j}\right\vert \max\limits_{1\leq j\leq n-1}\left\Vert \Delta
x_{j}\right\Vert ,\text{ \ \cite{DB}}; \\ 
\\ 
\dfrac{1}{2}\left( 1-\dfrac{1}{n}\right) \sum\limits_{j=1}^{n-1}\left\vert
\Delta \alpha _{j}\right\vert \sum\limits_{j=1}^{n-1}\left\Vert \Delta
x_{j}\right\Vert ,\text{ \ \ \ \ \ \cite{D1}}; \\ 
\\ 
\dfrac{1}{6}\cdot \dfrac{n^{2}-1}{n}\left( \sum\limits_{k=1}^{n-1}\left\vert
\Delta \alpha _{k}\right\vert ^{p}\right) ^{\frac{1}{p}}\left(
\sum\limits_{k=1}^{n-1}\left\Vert \Delta x_{k}\right\Vert ^{q}\right) ^{%
\frac{1}{q}} \\ 
\hfill \text{if }p>1,\ \frac{1}{p}+\frac{1}{q}=1,\text{ \ \ \ \cite{D2}.}%
\end{array}%
\right.   \notag
\end{align}%
The constants $\frac{1}{12},$ $\frac{1}{2}$ and $\frac{1}{6}$ are best
possible in the sense that they cannot be replaced by smaller constants.
\end{corollary}

In this paper, some new inequalities of Gr\"{u}ss type for sequences of
vectors in normed linear spaces subject of some boundedness conditions are
provided. Applications for discrete Fourier and Mellin transforms and for
vector-valued polynomials are pointed out as well.

\section{Some Analytic Inequalities}

The following result holds.

\begin{theorem}
\label{t2.1}Let $\left( X,\left\Vert \cdot \right\Vert \right) $ be a normed
linear space over the real or complex number field $\mathbb{K}$ $\left( 
\mathbb{K}=\mathbb{R},\mathbb{C}\right) ,$ $\overline{\mathbf{\alpha }}%
=\left( \alpha _{1},\dots ,\alpha _{n}\right) \in \mathbb{R}^{n},$ $%
\overline{\mathbf{p}}=\left( p_{1},\dots ,p_{n}\right) \in \mathbb{R}%
_{+}^{n} $ with $\sum_{i=1}^{n}p_{i}=1$ and $\overline{\mathbf{x}}=\left(
x_{1},\dots ,x_{n}\right) \in X^{n}.$

If $\alpha _{i}\in \overline{D}\left( \alpha ,R\right) :=\left\{ z\in 
\mathbb{K}|\left\vert z-\alpha \right\vert \leq R\right\} $ for some $\alpha
\in \mathbb{K}$ and $i\in \left\{ 1,\dots ,n\right\} ,$ $R>0,$ then we have
the inequality:%
\begin{equation}
\left\Vert \sum_{i=1}^{n}p_{i}\alpha _{i}x_{i}-\sum_{i=1}^{n}p_{i}\alpha
_{i}\cdot \sum_{i=1}^{n}p_{i}x_{i}\right\Vert \leq
R\sum_{i=1}^{n}p_{i}\left\Vert x_{i}-\sum_{j=1}^{n}p_{j}x_{j}\right\Vert .
\label{2.1}
\end{equation}%
The constant $c=1$ in the right hand side of the inequality is sharp in the
sense that it cannot be replaced by a smaller constant.
\end{theorem}

\begin{proof}
It is easy to see that, the following identity holds true%
\begin{equation}
\sum_{i=1}^{n}p_{i}\alpha _{i}x_{i}-\sum_{i=1}^{n}p_{i}\alpha _{i}\cdot
\sum_{i=1}^{n}p_{i}x_{i}=\sum_{i=1}^{n}p_{i}\left( \alpha _{i}-\alpha
\right) \left( x_{i}-\sum_{j=1}^{n}p_{j}x_{j}\right) .  \label{2.2}
\end{equation}%
Taking the norm in (\ref{2.2}), using the generalised triangle inequality
and the fact that $\alpha _{i}\in \overline{D}\left( \alpha ,R\right) ,$ $%
i=1,\dots ,n;$ we deduce%
\begin{align}
\left\Vert \sum_{i=1}^{n}p_{i}\alpha _{i}x_{i}-\sum_{i=1}^{n}p_{i}\alpha
_{i}\cdot \sum_{i=1}^{n}p_{i}x_{i}\right\Vert & \leq
\sum_{i=1}^{n}p_{i}\left\vert \alpha _{i}-\alpha \right\vert \left\Vert
x_{i}-\sum_{j=1}^{n}p_{j}x_{j}\right\Vert   \label{2.3} \\
& \leq R\sum_{i=1}^{n}p_{i}\left\Vert
x_{i}-\sum_{j=1}^{n}p_{j}x_{j}\right\Vert ,  \notag
\end{align}%
and the inequality (\ref{2.1}) is proved.

Now, assume that (\ref{2.1}) holds with a constant $C>0,$ i.e.,%
\begin{equation}
\left\Vert \sum_{i=1}^{n}p_{i}\alpha _{i}x_{i}-\sum_{i=1}^{n}p_{i}\alpha
_{i}\cdot \sum_{i=1}^{n}p_{i}x_{i}\right\Vert \leq
CR\sum_{i=1}^{n}p_{i}\left\Vert x_{i}-\sum_{j=1}^{n}p_{j}x_{j}\right\Vert ,
\label{2.4}
\end{equation}%
for $\overline{\mathbf{\alpha }}$, $\overline{\mathbf{p}}$ and $\overline{%
\mathbf{x}}$ as in the hypothesis of the theorem.

For $n=2,$ we have%
\begin{equation*}
\sum_{i=1}^{2}p_{i}\alpha _{i}x_{i}-\sum_{i=1}^{2}p_{i}\alpha _{i}\cdot
\sum_{i=1}^{2}p_{i}x_{i}=p_{2}p_{1}\left( \alpha _{2}-\alpha _{1}\right)
\left( x_{2}-x_{1}\right) 
\end{equation*}%
and%
\begin{equation*}
\sum_{i=1}^{2}p_{i}\left\Vert x_{i}-\sum_{j=1}^{2}p_{j}x_{j}\right\Vert
=2p_{2}p_{1}\left\Vert x_{2}-x_{1}\right\Vert ,
\end{equation*}%
and thus, by (\ref{2.4}), we deduce%
\begin{equation}
p_{2}p_{1}\left\vert \alpha _{2}-\alpha _{1}\right\vert \left\Vert
x_{2}-x_{1}\right\Vert \leq 2CRp_{2}p_{1}\left\Vert x_{2}-x_{1}\right\Vert .
\label{2.5}
\end{equation}%
If we choose $p_{2},p_{1}>0,$ $x_{1}\neq x_{2},$ $\alpha _{1}=\alpha -R,$ $%
\alpha _{2}=\alpha +R\in \overline{D}\left( \alpha ,R\right) ,$ then by (\ref%
{2.5}) we deduce $C\geq 1$ showing that $c_{0}=1$ is the best possible
constant in (\ref{2.1}).
\end{proof}

The following lemma holds.

\begin{lemma}
\label{l2.2}For the complex numbers $z,a,A\in \mathbb{C}$, the following
statements are equivalent

\begin{enumerate}
\item[(i)] $\func{Re}\left[ \left( A-z\right) \left( \overline{z}-\overline{a%
}\right) \right] \geq 0;$

\item[(ii)] $\left\vert z-\frac{a+A}{2}\right\vert \leq \frac{1}{2}%
\left\vert A-a\right\vert .$
\end{enumerate}
\end{lemma}

\begin{proof}
Define%
\begin{equation*}
I_{1}:=\func{Re}\left[ \left( A-z\right) \left( \overline{z}-\overline{a}%
\right) \right] =-\func{Re}\left( A\overline{a}\right) -\left\vert
z\right\vert ^{2}+\func{Re}\left[ z\overline{a}+\overline{z}A\right] 
\end{equation*}%
and%
\begin{align*}
I_{2}& :=\frac{1}{4}\left\vert A-a\right\vert ^{2}-\left\vert z-\frac{a+A}{2}%
\right\vert ^{2} \\
& =\frac{\left\vert A\right\vert ^{2}-2\func{Re}\left( A\overline{a}\right)
+\left\vert a\right\vert ^{2}}{4}-\left( \left\vert z\right\vert ^{2}-\func{%
Re}\left[ z\left( \overline{a}+\overline{A}\right) \right] +\frac{\left\vert
A+a\right\vert ^{2}}{4}\right)  \\
& =-\func{Re}\left[ A\overline{a}\right] -\left\vert z\right\vert ^{2}+\func{%
Re}\left( z\overline{a}\right) +\func{Re}\left( z\overline{A}\right)  \\
& =-\func{Re}\left[ A\overline{a}\right] -\left\vert z\right\vert ^{2}+\func{%
Re}\left( z\overline{a}\right) +\func{Re}\left( \overline{z}A\right) ,
\end{align*}%
since, obviously $\func{Re}\left( z\overline{A}\right) =\func{Re}\left( 
\overline{z\overline{A}}\right) .$

Consequently,%
\begin{equation}
\func{Re}\left[ \left( A-z\right) \left( \overline{z}-\overline{a}\right) %
\right] =\frac{1}{4}\left\vert A-a\right\vert ^{2}-\left\vert z-\frac{a+A}{2}%
\right\vert ^{2},  \label{2.6}
\end{equation}%
and the lemma is thus proved.
\end{proof}

\begin{remark}
\label{c2.3}For the real numbers $z,a,A\in \mathbb{R}$ (with $A\geq a$), the
following statements are obviously equivalent:
\end{remark}

\begin{corollary}
\begin{enumerate}
\item[(i)] $a\leq z\leq A;$

\item[(ii)] $\left\vert z-\frac{a+A}{2}\right\vert \leq \frac{A-a}{2}.$
\end{enumerate}
\end{corollary}

The following result of Gr\"{u}ss type for vectors in complex normed linear
spaces holds.

\begin{theorem}
\label{t2.4}Let $\left( X,\left\Vert \cdot \right\Vert \right) $ be a normed
linear space over the complex number field $\mathbb{C},$ $\overline{\mathbf{%
\alpha }}=\left( \alpha _{1},\dots ,\alpha _{n}\right) \in \mathbb{C}^{n},$ $%
\overline{\mathbf{p}}=\left( p_{1},\dots ,p_{n}\right) \in \mathbb{R}%
_{+}^{n} $ with $\sum_{i=1}^{n}p_{i}=1$ and $\overline{\mathbf{x}}=\left(
x_{1},\dots ,x_{n}\right) \in X^{n}.$

If there exists the complex numbers $a,A\in \mathbb{C}$ such that%
\begin{equation}
\func{Re}\left[ \left( A-\alpha _{i}\right) \left( \overline{\alpha _{i}}-%
\overline{a}\right) \right] \geq 0\text{ \ for each \ }i\in \left\{ 1,\dots
,n\right\}   \label{2.7}
\end{equation}%
or, equivalently,%
\begin{equation}
\left\vert \alpha _{i}-\frac{a+A}{2}\right\vert \leq \frac{1}{2}\left\vert
A-a\right\vert \text{ \ for each \ }i\in \left\{ 1,\dots ,n\right\} ,
\label{2.8}
\end{equation}%
then one has the inequality:%
\begin{equation}
\left\Vert \sum_{i=1}^{n}p_{i}\alpha _{i}x_{i}-\sum_{i=1}^{n}p_{i}\alpha
_{i}\cdot \sum_{i=1}^{n}p_{i}x_{i}\right\Vert \leq \frac{1}{2}\left\vert
A-a\right\vert \sum_{i=1}^{n}p_{i}\left\Vert
x_{i}-\sum_{j=1}^{n}p_{j}x_{j}\right\Vert .  \label{2.9}
\end{equation}%
The constant $\frac{1}{2}$ in the right hand side of the inequality is best
possible in the sense that it cannot be replaced by a smaller constant.
\end{theorem}

\begin{proof}
Follows by Theorem \ref{t2.1} on choosing $\alpha =\frac{a+A}{2}$ and $R=%
\frac{1}{2}\left\vert A-a\right\vert .$ The best constant may be shown in a
similar way as in the proof of Theorem \ref{t2.1}. We omit the details.
\end{proof}

The case of real normed linear spaces is embodied in the following corollary.

\begin{corollary}
\label{c2.5}Let $\left( X,\left\Vert \cdot \right\Vert \right) $ be a normed
linear space over the real number field $\mathbb{R},$ $\overline{\mathbf{%
\alpha }}=\left( \alpha _{1},\dots ,\alpha _{n}\right) \in \mathbb{R}^{n},$ $%
\overline{\mathbf{p}}=\left( p_{1},\dots ,p_{n}\right) \in \mathbb{R}%
_{+}^{n} $ with $\sum_{i=1}^{n}p_{i}=1$ and $\overline{\mathbf{x}}=\left(
x_{1},\dots ,x_{n}\right) \in X^{n}.$

If there exists the real numbers $m\leq M$ such that%
\begin{equation}
-\infty <m\leq a_{i}\leq M<\infty \text{ \ for each \ }i\in \left\{ 1,\dots
,n\right\} ,  \label{2.10}
\end{equation}%
then one has the inequality%
\begin{equation}
\left\Vert \sum_{i=1}^{n}p_{i}\alpha _{i}x_{i}-\sum_{i=1}^{n}p_{i}\alpha
_{i}\cdot \sum_{i=1}^{n}p_{i}x_{i}\right\Vert \leq \frac{1}{2}\left(
M-m\right) \sum_{i=1}^{n}p_{i}\left\Vert
x_{i}-\sum_{j=1}^{n}p_{j}x_{j}\right\Vert .  \label{2.11}
\end{equation}%
The constant $\frac{1}{2}$ is best possible in the sense mentioned above.
\end{corollary}

\begin{remark}
If $X=\mathbb{R}$, $\left\Vert .\right\Vert =\left\vert .\right\vert ,$ then
from $\left( \ref{2.11}\right) $ we obtain the inequality for real numbers
established for the first time in \cite{CD}.
\end{remark}

The dual result where some boundedness conditions for the sequence of
vectors are known, also holds.

\begin{theorem}
\label{t2.6}Let $X,$ $\overline{\mathbf{\alpha }}$, $\overline{\mathbf{p}}$
and $\overline{\mathbf{x}}$ be as in Theorem \ref{t2.1}.

If $x_{i}\in \overline{B}\left( x,R\right) :=\left\{ y\in X|\left\Vert
y-x\right\Vert \leq R\right\} $ for some $x\in X$ and $i\in \left\{ 1,\dots
,n\right\} ,$ $R>0,$ then we have the inequality:%
\begin{equation}
\left\Vert \sum_{i=1}^{n}p_{i}\alpha _{i}x_{i}-\sum_{i=1}^{n}p_{i}\alpha
_{i}\cdot \sum_{i=1}^{n}p_{i}x_{i}\right\Vert \leq
R\sum_{i=1}^{n}p_{i}\left\vert \alpha _{i}-\sum_{j=1}^{n}p_{j}\alpha
_{j}\right\vert .  \label{2.12}
\end{equation}%
The constant $c=1$ is sharp in the sense mentioned above.
\end{theorem}

\begin{proof}
It follows in a similar manner to the one in Theorem \ref{t2.1} on using the
following identity%
\begin{equation}
\sum_{i=1}^{n}p_{i}\alpha _{i}x_{i}-\sum_{i=1}^{n}p_{i}\alpha _{i}\cdot
\sum_{i=1}^{n}p_{i}x_{i}=\sum_{i=1}^{n}p_{i}\left( \alpha
_{i}-\sum_{j=1}^{n}p_{j}\alpha _{j}\right) \left( x_{i}-x\right) .
\label{2.13}
\end{equation}%
We omit the details.
\end{proof}

\begin{remark}
\label{r2.7}Using the Buniakowsky-Schwarz inequality for real numbers, we
may state that%
\begin{align*}
\sum_{i=1}^{n}p_{i}\left\vert \alpha _{i}-\sum_{j=1}^{n}p_{j}\alpha
_{j}\right\vert & \leq \left[ \sum_{i=1}^{n}p_{i}\left\vert \alpha
_{i}-\sum_{j=1}^{n}p_{j}\alpha _{j}\right\vert ^{2}\right] ^{\frac{1}{2}} \\
& =\left[ \sum_{i=1}^{n}p_{i}\left\vert \alpha _{i}\right\vert
^{2}-\left\vert \sum_{i=1}^{n}p_{i}\alpha _{i}\right\vert ^{2}\right] ^{%
\frac{1}{2}},
\end{align*}%
and then, by (\ref{2.12}) we may deduce the coarser bound%
\begin{align}
\left\Vert \sum_{i=1}^{n}p_{i}\alpha _{i}x_{i}-\sum_{i=1}^{n}p_{i}\alpha
_{i}\cdot \sum_{i=1}^{n}p_{i}x_{i}\right\Vert & \leq
R\sum_{i=1}^{n}p_{i}\left\vert \alpha _{i}-\sum_{j=1}^{n}p_{j}\alpha
_{j}\right\vert   \label{2.14} \\
& \leq R\left[ \sum_{i=1}^{n}p_{i}\left\vert \alpha _{i}\right\vert
^{2}-\left\vert \sum_{i=1}^{n}p_{i}\alpha _{i}\right\vert ^{2}\right] ^{%
\frac{1}{2}}.  \notag
\end{align}
\end{remark}

The following inequality for complex numbers, which is also interesting in
itself, holds.

\begin{proposition}
\label{t2.8}Let $\overline{\mathbf{\alpha }}=\left( \alpha _{1},\dots
,\alpha _{n}\right) \in \mathbb{C}^{n},$ $\overline{\mathbf{p}}=\left(
p_{1},\dots ,p_{n}\right) \in \mathbb{R}_{+}^{n}$ with $\sum_{i=1}^{n}p_{i}=1
$.
\end{proposition}

\begin{theorem}
If there exists the complex numbers $a,A\in \mathbb{C}$ such that (\ref{2.4}%
), or, equivalently (\ref{2.5}) holds, then one has the inequality%
\begin{equation}
0\leq \sum_{i=1}^{n}p_{i}\left\vert \alpha _{i}\right\vert ^{2}-\left\vert
\sum_{i=1}^{n}p_{i}\alpha _{i}\right\vert ^{2}\leq \frac{1}{4}\left\vert
A-a\right\vert ^{2}.  \label{2.15}
\end{equation}%
The constant $\frac{1}{4}$ is best possible in the sense that it cannot be
replaced by a smaler constant.
\end{theorem}

\begin{proof}
We apply Theorem \ref{t2.4} for the choice $X=\mathbb{C}$, $\left\Vert \cdot
\right\Vert =\left\vert \cdot \right\vert $ and $x_{i}=\overline{\alpha _{i}}
$ $\left( i=1,\dots ,n\right) .$ Then we get%
\begin{align*}
0& \leq \left[ \sum_{i=1}^{n}p_{i}\left\vert \alpha _{i}\right\vert
^{2}-\left\vert \sum_{i=1}^{n}p_{i}\alpha _{i}\right\vert ^{2}\right]  \\
& \leq \frac{1}{2}\left\vert A-a\right\vert \sum_{i=1}^{n}p_{i}\left\vert 
\overline{\alpha _{i}}-\sum_{j=1}^{n}p_{j}\overline{\alpha _{j}}\right\vert 
\\
& =\frac{1}{2}\left\vert A-a\right\vert \sum_{i=1}^{n}p_{i}\left\vert \alpha
_{i}-\sum_{j=1}^{n}p_{j}\alpha _{j}\right\vert  \\
& \leq \frac{1}{2}\left\vert A-a\right\vert \left[ \sum_{i=1}^{n}p_{i}\left%
\vert \alpha _{i}\right\vert ^{2}-\left\vert \sum_{i=1}^{n}p_{i}\alpha
_{i}\right\vert ^{2}\right] ^{\frac{1}{2}}
\end{align*}%
giving the desired result.

The fact that $\frac{1}{4}$ is the best possible constant may be proved in a
similar manner to the one incorporated in the proof of Theorem \ref{t2.1}.
\end{proof}

Another similar result for complex numbers also golds.

\begin{proposition}
\label{t2.9}With the assumptions of Proposition \ref{t2.8} for the complex
sequence $\overline{\mathbf{\alpha }}$, one has the inequality 
\begin{equation}
0\leq \left\vert \sum_{i=1}^{n}p_{i}\alpha _{i}^{2}-\left(
\sum_{i=1}^{n}p_{i}\alpha _{i}\right) ^{2}\right\vert \leq \frac{1}{4}%
\left\vert A-a\right\vert ^{2}.  \label{2.16}
\end{equation}%
The constant $\frac{1}{4}$ is best possible.
\end{proposition}

The proof follows by Theorem \ref{t2.4} and Remark \ref{r2.7} on choosing $X=%
\mathbb{C}$, $\left\Vert \cdot \right\Vert =\left\vert \cdot \right\vert $
and $x_{i}=\alpha _{i}$ $\left( i=1,\dots ,n\right) .$

\begin{remark}
\label{r2.10}Using the above results, we may state the following sequence of
inequalities of Gr\"{u}ss type for sequences of complex numbers%
\begin{align}
0& \leq \left\vert \sum_{i=1}^{n}p_{i}\alpha _{i}\beta
_{i}-\sum_{i=1}^{n}p_{i}\alpha _{i}\cdot \sum_{i=1}^{n}p_{i}\beta
_{i}\right\vert  \label{2.17} \\
& \leq \frac{1}{2}\left\vert A-a\right\vert \sum_{i=1}^{n}p_{i}\left\vert
\beta _{i}-\sum_{j=1}^{n}p_{j}\beta _{j}\right\vert  \notag \\
\text{(provided }\left\vert \alpha _{i}-\frac{a+A}{2}\right\vert & \leq 
\frac{1}{2}\left\vert A-a\right\vert \text{ for each }i\in \left\{ 1,\dots
,n\right\} \text{)}  \notag \\
& \leq \frac{1}{2}\left\vert A-a\right\vert \left(
\sum_{i=1}^{n}p_{i}\left\vert \beta _{i}\right\vert ^{2}-\left\vert
\sum_{i=1}^{n}p_{i}\beta _{i}\right\vert ^{2}\right) ^{\frac{1}{2}}  \notag
\\
& \leq \frac{1}{4}\left\vert A-a\right\vert \left\vert B-b\right\vert  \notag
\\
\text{(provided }\left\vert \beta _{i}-\frac{b+B}{2}\right\vert & \leq \frac{%
1}{2}\left\vert B-b\right\vert \text{ for each }i\in \left\{ 1,\dots
,n\right\} \text{).}  \notag
\end{align}%
The constants $\frac{1}{2}$ and $\frac{1}{4}$ are best possible in (\ref%
{2.17}).
\end{remark}

\section{Application for Discrete Fourier Transforms}

Let $\left( X,\left\Vert \cdot \right\Vert \right) $ be a normed linear
space over $\mathbb{C}$ and let $\overline{\mathbf{x}}=\left( x_{1},\dots
,x_{n}\right) $ be a sequence of vectors in $X.$

For a given $w\in \mathbb{R}$, define the \textit{discrete Fourier transform 
\cite{D2}}%
\begin{equation}
\mathcal{F}_{\omega }\left( \overline{\mathbf{x}}\right) \left( m\right)
:=\sum_{k=1}^{n}\exp \left( 2\omega imk\right) \cdot x_{k},\ \ \ \ \
m=1,\dots ,n;  \label{3.1}
\end{equation}%
where $i^{2}=-1.$

The following approximation result for the Fourier transform (\ref{3.1})
holds.

\begin{theorem}
\label{t3.1}If $x_{i}\in \overline{B}\left( x,R\right) :=\left\{ y\in
X|\left\Vert y-x\right\Vert \leq R\right\} ,$ $i\in \left\{ 1,\dots
,n\right\} ,$ for some $x\in X$ and $R>0,$ then we have the inequality:%
\begin{multline}
\left\Vert \mathcal{F}_{\omega }\left( \overline{\mathbf{x}}\right) \left(
m\right) -\frac{\sin \left( \omega mn\right) }{\sin \left( \omega m\right) }%
\exp \left[ \omega \left( n+1\right) im\right] \cdot \frac{1}{n}%
\sum_{k=1}^{n}x_{k}\right\Vert   \label{3.2} \\
\leq R\sum_{k=1}^{n}\left\vert \exp \left( 2\omega imk\right) -\frac{1}{n}%
\cdot \frac{\sin \left( \omega mn\right) }{\sin \left( \omega m\right) }\exp %
\left[ \omega \left( n+1\right) im\right] \right\vert ,
\end{multline}%
for all $m\in \left\{ 1,\dots ,n\right\} $ and $w\in \mathbb{R}$, $w\neq 
\frac{\ell }{m}\pi $, $\ell \in \mathbb{Z}$.
\end{theorem}

\begin{proof}
From the inequality (\ref{2.12}) of Theorem \ref{t2.6}, we have the
inequality%
\begin{equation}
\left\Vert \sum_{i=1}^{n}\alpha _{k}x_{k}-\sum_{i=1}^{n}\alpha _{k}\cdot 
\frac{1}{n}\sum_{i=1}^{n}x_{k}\right\Vert \leq R\sum_{k=1}^{n}\left\vert
\alpha _{k}-\frac{1}{n}\sum_{p=1}^{n}\alpha _{p}\right\vert ,  \label{3.3}
\end{equation}%
for any $\alpha _{k}\in \mathbb{C}$ and $x_{k}\in \overline{B}\left(
x,R\right) ,$ $k=1,\dots ,n.$

We may choose in (\ref{3.3}) $\alpha _{k}=\exp \left( 2wimk\right) $ to
obtain%
\begin{multline}
\left\Vert \mathcal{F}_{\omega }\left( \overline{\mathbf{x}}\right) \left(
m\right) -\sum_{k=1}^{n}\exp \left( 2\omega imk\right) \cdot \frac{1}{n}%
\sum_{k=1}^{n}x_{k}\right\Vert  \label{3.4} \\
\leq R\sum_{k=1}^{n}\left\vert \exp \left( 2\omega imk\right) -\frac{1}{n}%
\sum_{p=1}^{n}\exp \left( 2\omega imp\right) \right\vert
\end{multline}%
for all $m\in \left\{ 1,\dots ,n\right\} .$

Since, see for example \cite[p. 164]{D2}, by simple calculation we have%
\begin{equation*}
\sum_{k=1}^{n}\exp \left( 2\omega imk\right) =\frac{\sin \left( \omega
mn\right) }{\sin \left( \omega m\right) }\exp \left[ \omega \left(
n+1\right) im\right] ,
\end{equation*}%
for $w\neq \frac{\ell }{m}\pi ,$ $\ell \in \mathbb{Z}$, then by (\ref{3.4})
we deduce the desired inequality (\ref{3.2}).
\end{proof}

The following corollary is obvious.

\begin{corollary}
\label{c3.2}If $\overline{\mathbf{x}}=\left( x_{1},\dots ,x_{n}\right) \in 
\mathbb{C}^{n}$ and there exists $x,X\in \mathbb{C}$ such that%
\begin{equation}
\func{Re}\left[ \left( X-x_{i}\right) \left( \overline{x_{i}}-\overline{x}%
\right) \right] \geq 0\text{ \ for }i\in \left\{ 1,\dots ,n\right\} 
\label{3.5}
\end{equation}%
or, equivalently,%
\begin{equation}
\left\vert x_{i}-\frac{x+X}{2}\right\vert \leq \frac{1}{2}\left\vert
X-x\right\vert \text{ \ for }i\in \left\{ 1,\dots ,n\right\} ,  \label{3.6}
\end{equation}%
then we have the inequality%
\begin{multline}
\left\vert \mathcal{F}_{\omega }\left( \overline{\mathbf{x}}\right) \left(
m\right) -\frac{\sin \left( \omega mn\right) }{\sin \left( \omega m\right) }%
\exp \left[ \omega \left( n+1\right) im\right] \cdot \frac{1}{n}%
\sum_{k=1}^{n}x_{k}\right\vert   \label{3.7} \\
\leq \frac{1}{2}\left\vert X-x\right\vert \sum_{k=1}^{n}\left\vert \exp
\left( 2\omega imk\right) -\frac{1}{n}\cdot \frac{\sin \left( \omega
mn\right) }{\sin \left( \omega m\right) }\exp \left[ \omega \left(
n+1\right) im\right] \right\vert ,
\end{multline}%
for each $m\in \left\{ 1,\dots ,n\right\} $ and $w\in \mathbb{R}$, $w\neq 
\frac{\ell }{m}\pi ,$ $\ell \in \mathbb{Z}$.
\end{corollary}

\begin{remark}
If $\overline{\mathbf{x}}\in \mathbb{R}^{n}$ and there exists $a,A\in 
\mathbb{R}$ such that $a\leq x_{i}\leq A$ for $i\in \left\{ 1,\dots
,n\right\} $ then%
\begin{multline}
\left\vert \mathcal{F}_{\omega }\left( \overline{\mathbf{x}}\right) \left(
m\right) -\frac{\sin \left( \omega mn\right) }{\sin \left( \omega m\right) }%
\exp \left[ \omega \left( n+1\right) im\right] \cdot \frac{1}{n}%
\sum_{k=1}^{n}x_{k}\right\vert   \label{3.8} \\
\leq \frac{1}{2}\left( A-a\right) \sum_{k=1}^{n}\left\vert \exp \left(
2\omega imk\right) -\frac{1}{n}\cdot \frac{\sin \left( \omega mn\right) }{%
\sin \left( \omega m\right) }\exp \left[ \omega \left( n+1\right) im\right]
\right\vert 
\end{multline}%
for each $m\in \left\{ 1,\dots ,n\right\} $ and $w\in \mathbb{R}$, $w\neq 
\frac{\ell }{m}\pi ,$ $\ell \in \mathbb{Z}$.
\end{remark}

\section{Application for the Discrete Mellin Transform}

Let $\left( X,\left\Vert \cdot \right\Vert \right) $ be a normed linear
space over $\mathbb{K}$ ($\mathbb{K=C}$ or $\mathbb{K=R}$) and let $%
\overline{\mathbf{x}}=\left( x_{1},\dots ,x_{n}\right) $ be a sequence of
vectors in $X.$

Define the Mellin transform \cite{D2}%
\begin{equation}
\mathcal{M}\left( \overline{\mathbf{x}}\right) \left( m\right)
:=\sum_{k=1}^{n}k^{m-1}x_{k},\ \ \ \ m=1,\dots ,n,  \label{4.1}
\end{equation}%
where the sequence $\overline{\mathbf{x}}\in X^{n}.$

The following approximation result holds.

\begin{theorem}
\label{t4.1}If $x_{i}\in \overline{B}\left( x,R\right) ,$ $i\in \left\{
1,\dots ,n\right\} $ for some $x\in X$ and $R>0,$ then we have the
inequality:%
\begin{equation}
\left\Vert \mathcal{M}\left( \overline{\mathbf{x}}\right) \left( m\right)
-S_{m-1}\left( n\right) \cdot \frac{1}{n}\sum_{k=1}^{n}x_{k}\right\Vert \leq
R\sum_{k=1}^{n}\left\vert k^{m-1}-\frac{1}{n}S_{m-1}\left( n\right)
\right\vert ,  \label{4.2}
\end{equation}%
where $S_{p}\left( n\right) ,$ $p\in \mathbb{R}$, $n\in \mathbb{N}$ is the $%
p-$powered sum of the first $n$ natural numbers, i.e.,%
\begin{equation*}
S_{p}\left( n\right) :=\sum_{k=1}^{n}k^{p}.
\end{equation*}
\end{theorem}

\begin{proof}
We apply the inequality (\ref{3.3}) for $\alpha _{k}=k^{m-1}$ to obtain%
\begin{equation}
\left\Vert \sum_{k=1}^{n}k^{m-1}x_{k}-\sum_{k=1}^{n}k^{m-1}\cdot \frac{1}{n}%
\sum_{k=1}^{n}x_{k}\right\Vert \leq R\sum_{k=1}^{n}\left\vert k^{m-1}-\frac{1%
}{n}\sum_{l=1}^{n}l^{m-1}\right\vert ,  \label{4.3}
\end{equation}%
giving the desired result (\ref{4.2}).
\end{proof}

For $m=2,$ we have 
\begin{align*}
\sum_{k=1}^{n}\left\vert k-\frac{1}{n}S_{1}\left( n\right) \right\vert &
=\sum_{k=1}^{n}\left\vert k-\frac{n+1}{2}\right\vert  \\
& =\sum_{k=1}^{\left[ \frac{n+1}{2}\right] }\left( \frac{n+1}{2}-k\right)
+\sum_{k=\left[ \frac{n+1}{2}\right] +1}^{k}\left( k-\frac{n+1}{2}\right)  \\
& =:I,
\end{align*}%
where $\left[ a\right] $ is the integer part of $a\in \mathbb{R}$.

Observe that%
\begin{equation*}
\sum_{k=1}^{\left[ \frac{n+1}{2}\right] }\left( \frac{n+1}{2}-k\right) =%
\frac{n+1}{2}\left[ \frac{n+1}{2}\right] -\frac{\left[ \frac{n+1}{2}\right]
\left( \left[ \frac{n+1}{2}\right] +1\right) }{2}
\end{equation*}%
and%
\begin{align*}
\sum_{k=\left[ \frac{n+1}{2}\right] +1}^{n}\left( k-\frac{n+1}{2}\right) &
=\sum_{k=1}^{n}k-\sum_{k=1}^{\left[ \frac{n+1}{2}\right] }k-\frac{n+1}{2}%
\left( n-\left[ \frac{n+1}{2}\right] \right)  \\
& =\frac{n+1}{2}\left[ \frac{n+1}{2}\right] -\frac{\left[ \frac{n+1}{2}%
\right] \left( \left[ \frac{n+1}{2}\right] +1\right) }{2},
\end{align*}%
thus%
\begin{equation*}
I=\left( n+1\right) \left[ \frac{n+1}{2}\right] -\left[ \frac{n+1}{2}\right]
\left( \left[ \frac{n+1}{2}\right] +1\right) =\left[ \frac{n+1}{2}\right]
\left( n-\left[ \frac{n+1}{2}\right] \right) .
\end{equation*}

Now, if we consider a particular value of the Mellin transform%
\begin{equation*}
\mu \left( \overline{\mathbf{x}}\right) :=\sum_{k=1}^{n}kx_{k},
\end{equation*}%
then we may state the following.

\begin{corollary}
\label{c4.2}With the assumptions of Theorem \ref{t4.1}, we have%
\begin{equation}
\left\Vert \mu \left( \overline{\mathbf{x}}\right) -\frac{n+1}{2}%
\sum_{k=1}^{n}x_{k}\right\Vert \leq R\cdot \left[ \frac{n+1}{2}\right]
\left( n-\left[ \frac{n+1}{2}\right] \right) .  \label{4.4}
\end{equation}
\end{corollary}

\begin{remark}
\label{r4.3}Assume that $\overline{\mathbf{x}}=\left( x_{1},\dots
,x_{n}\right) \in \mathbb{C}^{n}$ are such that there exists $x,X\in \mathbb{%
C}$ such that (\ref{3.5}) or, equivalently, (\ref{3.6}) holds. Then we have
the inequality%
\begin{equation}
\left\vert \mathcal{M}\left( \overline{\mathbf{x}}\right) \left( m\right)
-S_{m-1}\left( n\right) \cdot \frac{1}{n}\sum_{k=1}^{n}x_{k}\right\vert \leq 
\frac{1}{2}\left\vert X-x\right\vert \sum_{k=1}^{n}\left\vert k^{m-1}-\frac{1%
}{n}S_{m-1}\left( n\right) \right\vert .  \label{4.5}
\end{equation}%
In particular, we have the inequality%
\begin{equation}
\left\vert \mu \left( \overline{\mathbf{x}}\right) -\frac{n+1}{2}%
\sum_{k=1}^{n}x_{k}\right\vert \leq \frac{1}{2}\left\vert X-x\right\vert %
\left[ \frac{n+1}{2}\right] \left( n-\left[ \frac{n+1}{2}\right] \right) .
\label{4.6}
\end{equation}
\end{remark}

\begin{remark}
\label{r4.4}Assume that $\overline{\mathbf{x}}=\left( x_{1},\dots
,x_{n}\right) \in \mathbb{R}^{n}$ such that there exists $a,A\in \mathbb{R}$
with $a\leq x_{i}\leq A$, $i\in \left\{ 1,\dots ,n\right\} .$ Then we have
the inequality%
\begin{equation}
\left\vert \mathcal{M}\left( \overline{\mathbf{x}}\right) \left( m\right)
-S_{m-1}\left( n\right) \cdot \frac{1}{n}\sum_{k=1}^{n}x_{k}\right\vert \leq 
\frac{1}{2}\left( A-a\right) \sum_{k=1}^{n}\left\vert k^{m-1}-\frac{1}{n}%
S_{m-1}\left( n\right) \right\vert   \label{4.7}
\end{equation}%
in particular, we have the inequality%
\begin{equation}
\left\vert \mu \left( \overline{\mathbf{x}}\right) -\frac{n+1}{2}%
\sum_{k=1}^{n}x_{k}\right\vert \leq \frac{1}{2}\left( A-a\right) \left[ 
\frac{n+1}{2}\right] \left( n-\left[ \frac{n+1}{2}\right] \right) .
\label{4.8}
\end{equation}
\end{remark}

\section{Application for Polynomials}

Let $\left( X,\left\Vert \cdot \right\Vert \right) $ be a normed linear
space over $\mathbb{C}$ and let $\overline{\mathbf{c}}=\left( c_{0},\dots
,c_{n}\right) $ be a sequence of vectors in $X.$

Define the polynomial $P:\mathbb{C\rightarrow }X$ with the coefficients $%
\overline{\mathbf{c}}=\left( c_{0},\dots ,c_{n}\right) $ by%
\begin{equation*}
P\left( z\right) =c_{0}+zc_{1}+z^{2}c_{2}+\cdots +z^{n}c_{n},\ \ \ \ z\in 
\mathbb{C},\ \ c_{n}\neq 0.
\end{equation*}

The following approximation holds.

\begin{theorem}
\label{t5.1}If $c_{i}\in \overline{B}\left( x_{0},R\right) ,$ $i\in \left\{
0,\dots ,n\right\} $ for some $x_{0}\in X$ and $R>0,$ then we have the
inequality%
\begin{equation}
\left\Vert P\left( z\right) -\frac{z^{n+1}-1}{z-1}\cdot \frac{1}{n+1}%
\sum_{k=0}^{n}c_{k}\right\Vert \leq R\sum_{k=0}^{n}\left\vert z^{k}-\frac{1}{%
n+1}\cdot \frac{z^{n+1}-1}{z-1}\right\vert  \label{5.1}
\end{equation}%
for all $z\in \mathbb{C}$, $z\neq 1.$
\end{theorem}

\begin{proof}
The proof follows by the inequality (\ref{3.3}) on choosing $\alpha
_{k}=z^{k},$ and $x_{k}=c_{k},$ $k=0,\dots ,n$ and we omit the details.
\end{proof}

The following corollary concerning the location of $P\left( z_{k}\right) ,$
where $z_{k}$ are the complex roots of the unity holds.

\begin{corollary}
\label{c5.2}Let $z_{k}=\cos \left( \frac{k\pi }{n+1}\right) +i\sin \left( 
\frac{k\pi }{n+1}\right) ,$ where $k\in \left\{ 0,\dots ,n\right\} ,$ be the
complex $\left( n+1\right) -$roots of the unity. Then we have the inequality%
\begin{equation}
\left\Vert P\left( z_{n}\right) \right\Vert \leq \left( n+1\right) R,\ \ \ \
\ z\in \left\{ 1,\dots ,n\right\} ,  \label{5.2}
\end{equation}%
where the coefficients $c_{i}$ $\left( i\in \left\{ 0,\dots ,n\right\}
\right) $ satisfy the assumptions of Theorem \ref{t5.1}.
\end{corollary}

\begin{proof}
Follows by (\ref{5.1}) on choosing $z=z_{n},$ $k\in \left\{ 1,\dots
,n\right\} $ and taking into account that $z_{k}^{n+1}=1$ and $\left\vert
z_{k}\right\vert ^{k}=1$ for $k\in \left\{ 1,\dots ,n\right\} .$
\end{proof}

\begin{remark}
\label{r5.3}Assume that $\overline{\mathbf{c}}=\left( c_{0},\dots
,c_{n}\right) \in \mathbb{C}^{n+1}$ are such that there exists $w,W\in 
\mathbb{C}$ with the property%
\begin{equation}
\func{Re}\left[ \left( W-c_{i}\right) \left( \overline{c_{i}}-\overline{w}%
\right) \right] \geq 0\text{ \ for }i\in \left\{ 0,\dots ,n\right\}
\label{5.3}
\end{equation}%
or, equivalently,%
\begin{equation}
\left\vert c_{i}-\frac{w+W}{2}\right\vert \leq \frac{1}{2}\left\vert
W-w\right\vert \text{ \ for }i\in \left\{ 0,\dots ,n\right\} .  \label{5.4}
\end{equation}%
Then we have the inequality%
\begin{equation}
\left\vert P\left( z\right) -\frac{z^{n+1}-1}{z-1}\cdot \frac{1}{n+1}%
\sum_{k=0}^{n}c_{k}\right\vert \leq \frac{1}{2}\left\vert W-w\right\vert
\sum_{k=0}^{n}\left\vert z^{k}-\frac{1}{n+1}\cdot \frac{z^{n+1}-1}{z-1}%
\right\vert  \label{5.5}
\end{equation}%
for any $z\in \mathbb{C}$, $z\neq 1.$

If $z_{k},$ $k\in \left\{ 0,\dots ,n\right\} $ are the $\left( n+1\right) -$%
roots of the unity, then%
\begin{equation}
\left\vert P\left( z_{k}\right) \right\vert \leq \frac{1}{2}\left\vert
W-w\right\vert \left( n+1\right) ,\ \ \ \ \ z\in \left\{ 1,\dots ,n\right\} ,
\label{5.6}
\end{equation}%
provided the complex coefficients of $P\left( z\right) $ satisfy either (\ref%
{5.3}) or (\ref{5.4}).
\end{remark}

\begin{remark}
\label{r5.4}Assume that the coefficients \thinspace $c_{i}$ $\left( i\in
\left\{ 0,\dots ,n\right\} \right) $ are real numbers with the property that
there exists $a,A\in \mathbb{R}$ such that $a\leq c_{i}\leq A,$ $i\in
\left\{ 0,\dots ,n\right\} .$ Then we have the inequality%
\begin{equation}
\left\vert P\left( z\right) -\frac{z^{n+1}-1}{z-1}\cdot \frac{1}{n+1}%
\sum_{k=0}^{n}c_{k}\right\vert \leq \frac{1}{2}\left( A-a\right)
\sum_{k=0}^{n}\left\vert z^{k}-\frac{1}{n+1}\cdot \frac{z^{n+1}-1}{z-1}%
\right\vert  \label{5.7}
\end{equation}%
for any $z\in \mathbb{C}$, $z\neq 1.$

In particular, if $z_{k},$ $k\in \left\{ 0,\dots ,n\right\} $ are the $%
\left( n+1\right) -$roots of the unity and the coefficients of $P\left(
z\right) $ satisfy the above assumption, then%
\begin{equation}
\left\vert P\left( z_{k}\right) \right\vert \leq \frac{1}{2}\left(
A-a\right) \left( n+1\right) .  \label{5.8}
\end{equation}
\end{remark}


\begin{thebibliography}{9}
\bibitem{CD} P. CERONE and S.S. DRAGOMIR, A refinement of the Gr\"{u}ss
inequality and applications, \textit{RGMIA Res. Rep. Coll.}, \textbf{5}%
(2002), No. 2, Article 14. On line: \texttt{\
http://rgmia.vu.edu.au/v5n2.html}

\bibitem{DB} S.S. DRAGOMIR and\ G.L. BOOTH, Gr\"{u}ss-Lupa\c{s} type
inequality and its application for $p$-moments of guessing mappings, \textit{%
Math. Comm., }\textbf{5} (2000), 117-126.

\bibitem{D1} S.S. DRAGOMIR, A Gr\"{u}ss type inequality for sequences of
vectors in normed linear spaces and applications, \textit{RGMIA Res. Rep.
Coll}., \textbf{5}(2) (2002), Article 9. Available online at \texttt{%
http://rgmia.vu.edu.au/v5n2.html}

\bibitem{D2} S.S. DRAGOMIR, Another Gr\"{u}ss type inequality for sequences
of vectors in normed linear spaces and applications, \textit{J. Comp. Anal.
\& Appl., }\textbf{4}(2) (2002), 155-172.
\end{thebibliography}
\end{document}